# LOG-SOBOLEV INEQUALITIES: DIFFERENT ROLES OF RIC AND HESS

By Feng-Yu Wang[1]

*Beijing Normal University and Swansea University*

Let $P_t$ be the diffusion semigroup generated by $L := \Delta + \nabla V$ on a complete connected Riemannian manifold with $\mathrm{Ric} \geq -(\sigma^2 \rho_o^2 + c)$ for some constants $\sigma, c > 0$ and $\rho_o$ the Riemannian distance to a fixed point. It is shown that $P_t$ is hypercontractive, or the log-Sobolev inequality holds for the associated Dirichlet form, provided $-\mathrm{Hess}_V \geq \delta$ holds outside of a compact set for some constant $\delta > (1+\sqrt{2})\sigma\sqrt{d-1}$. This indicates, at least in finite dimensions, that Ric and $-\mathrm{Hess}_V$ play quite different roles for the log-Sobolev inequality to hold. The supercontractivity and the ultracontractivity are also studied.

**1. Introduction.** Let $M$ be a $d$-dimensional completed connected noncompact Riemannian manifold and $V \in C^2(M)$ such that

$$Z := \int_M e^{V(x)}\,dx < \infty, \tag{1.1}$$

where $dx$ is the volume measure on $M$. Let $\mu(dx) = Z^{-1} e^{V(x)}\,dx$. Under (1.1) it is easy to see that $H_0^{2,1}(\mu) = W^{2,1}(\mu)$, where $H_0^{2,1}(\mu)$ is the completion of $C_0^1(M)$ under the Sobolev norm $\|f\|_{2,1} := \mu(f^2 + |\nabla f|^2)^{1/2}$, and $W^{2,1}(\mu)$ is the completion of the class $\{f \in C^1(M): f + |\nabla f| \in L^2(\mu)\}$ under $\|\cdot\|_{2,1}$. Then the $L$-diffusion process is nonexplosive and its semigroup $P_t$ is uniquely determined. Moreover, $P_t$ is symmetric in $L^2(\mu)$ so that $\mu$ is $P_t$-invariant. It is well known by the Bakry–Emery criterion (see [4]) that

$$\mathrm{Ric} - \mathrm{Hess}_V \geq K \tag{1.2}$$

Received October 2007; revised August 2008.
[1]Supported in part by WIMCS, Creative Research Group Fund of the National Natural Science Foundation of China (No. 10721091) and the 973-Project.

for some constant $K > 0$ implies the Gross log-Sobolev inequality [14],

$$\mu(f^2 \log f^2) := \int_M f^2 \log f^2 \, d\mu \leq C\mu(|\nabla f|^2),$$
(1.3)
$$\mu(f^2) = 1, f \in C^1(M)$$

for $C = 2/K$. This result was extended by Chen and the author [9] to the situation that $\text{Ric} - \text{Hess}_V$ is uniformly positive outside a compact set. In the case that $\text{Ric} - \text{Hess}_V$ is bounded below, sufficient concentration conditions of $\mu$ for (1.3) to hold are presented in [1, 19, 20]. Obviously, in a condition on $\text{Ric} - \text{Hess}_V$ the Ricci curvature and $-\text{Hess}_V$ play the same role.

What can we do when $\text{Ric} - \text{Hess}_V$ is unbounded below? It seems very hard to confirm the log-Sobolev inequality with the unbounded below condition of $\text{Ric} - \text{Hess}_V$. Therefore, in this paper we try to clarify the roles of Ric and $-\text{Hess}_V$ in the study of the log-Sobolev inequality. Let us first recall the gradient estimate of $P_t$, which is a key point in the above references to prove the log-Sobolev inequality.

Let $x_t$ be the $L$-diffusion process starting at $x$, and let $v \in T_x M$. Due to Bismut [6] and Elworthy–Li [11], under a reasonable lower bound condition of $\text{Ric} - \text{Hess}_V$, one has

$$\langle \nabla P_t f, v \rangle = \mathbb{E} \langle \nabla f(x_t), v_t \rangle, \qquad t > 0, f \in C_b^1(M),$$

where $v_t \in T_{x_t} M$ solves the equation

$$D_t v_t := //_{t \to 0}^{-1} \frac{d}{dt} //_{t \to 0} v_t = -(\text{Ric} - \text{Hess}_V)^\#(v_t)$$

for $//_{t \to 0} : T_{x_t} M \to T_x M$ the associated stochastic parallel displacement, and $(\text{Ric} - \text{Hess}_V)^\#(v_t) \in T_{x_t} M$ with

$$\langle (\text{Ric} - \text{Hess}_V)^\#(v_t), X \rangle := (\text{Ric} - \text{Hess}_V)(v_t, X), \qquad X \in T_{x_t} M.$$

Thus, for the gradient of $P_t$, which is a short distance behavior of the diffusion process, a condition on $\text{Ric} - \text{Hess}_V$ appears naturally.

On the other hand, however, Ric and $-\text{Hess}_V$ play very different roles for long distance behaviors. For instance, Let $\rho_o$ be the Riemannian distance function to a fixed point $o \in M$. If $\text{Ric} \geq -k$ and $-\text{Hess}_V \geq \delta$ for some $k \geq 0$, $\delta \in \mathbb{R}$, the Laplacian comparison theorem implies

$$L\rho_o \leq \sqrt{k(d-1)} \coth[\sqrt{k/(d-1)}\rho_o] - \delta \rho_o.$$

Therefore, for large $\rho_o$, the Ric lower bound leads to a bounded term while that of $-\text{Hess}_V$ provides a linear term. The same phenomena appears in the formula on distance of coupling by parallel displacement (cf. [3], (2.3), (2.4)), which implies the above Bismut–Elworthy–Li formula by letting the initial distance tend to zero (cf. [15]). Here, $k \geq 0$ is essential for our framework, since the manifold has to be compact, if Ric is bounded below by a positive constant.

Since the log-Sobolev inequality is always available on bounded regular domains, it is more likely a long-distance property of the diffusion process. So, Ric and $-\text{Hess}_V$ should take different roles in the study of the log-Sobolev inequality. Indeed, it has been observed by the author [20] that (1.3) holds for some $C > 0$, provided Ric is bounded below and $-\text{Hess}_V$ is uniformly positive outside a compact set. This indicates that for the log-Sobolev inequality, the positivity of $-\text{Hess}_V$ is a dominative condition, which allows the Ricci curvature to be bounded below by an arbitrary negative constant, and hence, allows $\text{Ric} - \text{Hess}_V$ to be globally negative on $M$.

The first aim of this paper is to search for the weakest possibility of curvature lower bound for the log-Sobolev inequality to hold under the condition

(1.4) $\qquad -\text{Hess}_V \geq \delta \qquad$ outside a compact set

for some constant $\delta > 0$. This condition is reasonable as the log-Sobolev inequality implies $\mu(e^{\lambda \rho_o^2}) < \infty$ for some $\lambda > 0$ (see, e.g., [2, 17]).

According to the following Theorem 1.1 and Example 1.1, we conclude that under (1.4) the optimal curvature lower bound condition for (1.3) to hold is

(1.5) $$\inf_M \{\text{Ric} + \sigma^2 \rho_o^2\} > -\infty$$

for some constant $\sigma > 0$, such that $\delta > (1 + \sqrt{2})\sigma\sqrt{d-1}$. More precisely, let $\theta_0 > 0$ be the smallest positive constant, such that for any connected complete noncompact Riemannian manifold $M$ and $V \in C^2(M)$, such that $Z := \int_M e^{V(x)}\,dx < \infty$, the conditions (1.4) and (1.5) with $\delta > \sigma\theta_0\sqrt{d-1}$, implies (1.3) for some $C > 0$. Due to Theorem 1.1 and Example 1.1 below, we conclude that

$$\theta_0 \in [1, 1+\sqrt{2}].$$

The exact value of $\theta_0$ is however unknown.

THEOREM 1.1. *Assume that (1.4) and (1.5) hold for some constants $c, \delta, \sigma > 0$ with $\delta > (1+\sqrt{2})\sigma\sqrt{d-1}$. Then (1.3) holds for some $C > 0$.*

EXAMPLE 1.1. Let $M = \mathbb{R}^2$ be equipped with the rotationally symmetric metric

$$ds^2 = dr^2 + \{re^{kr^2}\}^2\,d\theta^2,$$

under the polar coordinates $(r, \theta) \in [0, \infty) \times \mathbb{S}^1$ at 0, where $k > 0$ is a constant, then (see, e.g., [13])

$$\text{Ric} = -\frac{(d^2/dr^2)(re^{kr^2})}{re^{kr^2}} = -4k - 4k^2r^2.$$

Thus, (1.5) holds for $\sigma = 2k$. Next, take $V = -k\rho_o^2 - \lambda(\rho_o^2 + 1)^{1/2}$ for some $\lambda > 0$. By the Hessian comparison theorem and the negativity of the sectional curvature, we obtain (1.4) for $\delta = 2k$. Since $d = 2$ and

$$(1.6) \qquad e^{V(x)}\,dx = re^{-\lambda(1+r^2)^{1/2}}\,dr\,d\theta,$$

one has $Z < \infty$ and $\delta = 2k = \sigma\sqrt{d-1}$. But the log-Sobolev inequality is not valid since by Herbst's inequality it implies $\mu(e^{r\rho_o^2}) < \infty$ for some $r > 0$, which is, however, not the case due to (1.6). Since in this example one has $\delta > \sigma\theta\sqrt{d-1}$ for any $\theta < 1$, according to the definition of $\theta_0$, we conclude that $\theta_0 \geq 1$.

Following the line of [19, 20], the key point in the proof of Theorem 1.1 will be a proper Harnack inequality of type

$$(P_t f(x))^\alpha \leq C_\alpha(t,x,y) P_t f^\alpha(y), \qquad t > 0, x, y \in M,$$

for any nonnegative $f \in C_b(M)$, where $\alpha > 1$ is a constant and $C_\alpha \in C((0,\infty), M^2)$ is a positive function. Such an inequality was established in [19] for $\mathrm{Ric} - \mathrm{Hess}_V$ bounded below and extended in [3] to a more general situation with Ric satisfying (1.5).

The Harnack inequality presented in [3] contains a leading term $\exp[\rho(x,y)^4]$, which is, however, too large to be integrability w.r.t. $\mu \times \mu$ under our conditions. So, to prove Theorem 1.1, we shall present a sharper Harnack inequality in Section 3 by refining the coupling method introduced in [3] (see Proposition 3.1 below). This inequality, together with the concentration of $\mu$ ensured by (1.4) and (1.5), will imply the hypercontractivity of $P_t$. To establish this new Harnack inequality, some necessary preparations are presented in Section 2.

Finally, in the same spirit of Theorem 1.1, the supercontractivity and ultracontractivity of $P_t$ are studied in Section 4 under explicit conditions on Ric and $-\mathrm{Hess}_V$.

**2. Preparations.** We first study the concentration of $\mu$ by using (1.4) and (1.5), for which we need to estimate $L\rho_o$ from above according to [5] and references within.

LEMMA 2.1. *If (1.4) and (1.5) hold, then there exists a constant $C_1 > 0$ such that*

$$(2.1) \qquad L\rho_o^2 \leq C_1(1 + \rho_o) - 2(\delta - \sigma\sqrt{d-1})\rho_o^2$$

*holds outside* $\mathrm{cut}(o)$, *the cut-locus of $o$. If moreover $\delta > \sigma\sqrt{d-1}$ then $Z < \infty$ and $\mu(e^{\lambda\rho_o^2}) < \infty$ for all $\lambda < \frac{1}{2}(\delta - \sigma\sqrt{d-1})$.*

PROOF. By (1.5) we have $\mathrm{Ric} \geq -(c+\sigma^2\rho_o^2)$ for some constant $c>0$. By the Laplacian comparison theorem this implies that

$$\Delta\rho_o \leq \sqrt{(c+\sigma^2\rho_o^2)(d-1)}\coth[\sqrt{(c+\sigma^2\rho_o^2)/(d-1)}\ \rho_o]$$

holds outside $\mathrm{cut}(o)$. Thus, outside $\mathrm{cut}(o)$ one has

$$\begin{aligned}\Delta\rho_o^2 &\leq 2\rho_o\sqrt{(c+\sigma^2\rho_o^2)(d-1)}\coth[\sqrt{(c+\sigma^2\rho_o^2)/(d-1)}\ \rho_o]+2\\ &\leq 2d+2\rho_o\sqrt{(c+\sigma^2\rho_o^2)(d-1)},\end{aligned}$$
(2.2)

where the second inequality follows from the fact that

$$r\cosh r \leq (1+r)\sinh r, \qquad r\geq 0.$$

On the other hand, for $x \notin \mathrm{cut}(o)$ and $U$ the unit tangent vector along the unique minimal geodesic $\ell$ form $o$ to $x$, by (1.4) there exists a constant $c_1>0$ independent of $x$ such that

$$\langle\nabla V,\nabla\rho_o\rangle(x) = \langle\nabla V,U\rangle(o)+\int_0^{\rho_o(x)}\mathrm{Hess}_V(U,U)(\ell_s)\,ds \leq c_1-\delta\rho_o(x).$$

Combining this with (2.2) we prove (2.1).

Finally, let $\delta > \sigma\sqrt{d-1}$ and $0<\lambda<\frac{1}{2}(\delta-\sigma\sqrt{d-1})$. By (2.1) we have

$$\begin{aligned}Le^{\lambda\rho_o^2} &\leq \lambda e^{\lambda\rho_o^2}(C_1(1+\rho_o)-2(\delta-\sigma\sqrt{d-1})\rho_o^2+4\lambda\rho_o^2)\\ &\leq c_2-c_3\rho_o^2 e^{\lambda\rho_o^2}\end{aligned}$$

for some constants $c_2,c_3>0$. By [5], Proposition 3.2, this implies $Z<\infty$ and

$$\int_M \rho_o^2 e^{\lambda\rho_o^2}\,d\mu \leq \frac{c_2}{c_3}<\infty. \qquad \square$$

LEMMA 2.2. *Let $x_t$ be the $L$-diffusion process with $x_0=x\in M$. If (1.4) and (1.5) hold with $\delta>\sigma\sqrt{d-1}$, then for any $\delta_0\in(\sigma\sqrt{d-1},\delta)$ there exists a constant $C_2>0$ such that*

$$\begin{aligned}\mathbb{E}\exp&\left[\frac{(\delta_0-\sigma\sqrt{d-1})^2}{4}\int_0^T\rho_o(x_t)^2\,dt\right]\\ &\leq \exp\left[C_2T+\frac{1}{4}(\delta_0-\sigma\sqrt{d-1})\rho_o(x)^2\right], \qquad T>0, x\in M.\end{aligned}$$

PROOF. By Lemma 2.1, we have

$$L\rho_o^2 \leq C-2(\delta_0-\sigma\sqrt{d-1})\rho_o^2$$

outside cut($o$) for some constant $C > 0$. Then the Itô formula for $\rho_o(x_t)$ due to Kendall [16] implies that

$$(2.3) \quad d\rho_o^2(x_t) \leq 2\sqrt{2}\rho_o(x_t)\,db_t + [C - 2(\delta_0 - \sigma\sqrt{d-1})\rho_o^2(x_t)]\,dt$$

holds for some Brownian motion $b_t$ on $\mathbb{R}$. This implies that the $L$-diffusion process is nonexplosive so that

$$T_n := \inf\{t \geq 0 : \rho_o(x_t) \geq n\} \to \infty$$

as $n \to \infty$. Indeed, (2.3) implies that

$$n\mathbb{P}(T_n \leq t) \leq \mathbb{E}\rho_o(x_{t \wedge T_n})^2 \leq \rho_o(x)^2 + Ct, \qquad n \geq 1, t > 0.$$

Hence, $\mathbb{P}(T_n \leq t) \to 0$ as $n \to \infty$ for any $t > 0$. This implies $\lim_{n \to \infty} T_n = \infty$ a.s.

For any $\lambda > 0$ and $n \geq 1$, it follows from (2.3) that

$$\mathbb{E}\exp\left[2\lambda(\delta_0 - \sigma\sqrt{d-1})\int_0^{T \wedge T_n} \rho_o^2(x_t)\,dt\right]$$

$$\leq e^{\lambda\rho_o^2(x) + C\lambda T}\mathbb{E}\exp\left[2\sqrt{2}\lambda\int_0^{T \wedge T_n}\rho_o(x_t)\,db_t\right]$$

$$\leq e^{\lambda\rho_o^2(x) + C\lambda T}\left(\mathbb{E}\exp\left[16\lambda^2\int_0^{T \wedge T_n}\rho_o^2(x_t)\,dt\right]\right)^{1/2},$$

where in the last step we have used the inequality

$$\mathbb{E}e^{M_t} \leq (\mathbb{E}e^{2\langle M\rangle_t})^{1/2}$$

for $M_t = 2\sqrt{2}\lambda \int_0^{t \wedge T_n}\rho_o(X_s)\,db_s$. This follows immediately from the Schwartz inequality and the fact that $\exp[2M_t - 2\langle M\rangle_t]$ is a martingale. Thus, taking

$$\lambda = \tfrac{1}{8}(\delta_0 - \sigma\sqrt{d-1}),$$

we obtain

$$\mathbb{E}\exp\left[\frac{1}{4}(\delta_0 - \sigma\sqrt{d-1})^2\int_0^{T \wedge T_n}\rho_o^2(x_t)\,dt\right]$$

$$\leq \exp\left[\frac{1}{4}(\delta_0 - \sigma\sqrt{d-1})\rho_o^2(x) + C_2 T\right]$$

for some $C_2 > 0$. Then the proof is completed by letting $n \to \infty$. □

Finally, we recall the coupling argument introduced in [3] for establishing the Harnack inequality of $P_t$.

Let $T > 0$ and $x \neq y \in M$ be fixed. Then the $L$-diffusion process starting from $x$ can be constructed by solving the following Itô stochastic differential equation:

$$d_I x_t = \sqrt{2}\Phi_t\,dB_t + \nabla V(x_t)\,dt, \qquad x_0 = x,$$

where $d_I$ is the Itô differential on manifolds introduced in [12] (see also [3]), $B_t$ is the $d$-dimensional Brownian motion, and $\Phi_t$ is the horizontal lift of $x_t$ onto the orthonormal frame bundle $O(M)$.

To construct another diffusion process $y_t$ starting from $y$ such that $x_T = y_T$, as in [3], we add an additional drift term to the equation (as explained in [3], Section 3, we may and do assume that the cut-locus of $M$ is empty)

$$d_I y_t = \sqrt{2} P_{x_t,y_t} \Phi_t \, dB_t + \nabla V(y_t) \, dt + \xi_t U(x_t, y_t) 1_{\{t<\tau\}} \, dt, \qquad y_0 = y,$$

where $P_{x_t,y_t}$ is the parallel transformation along the unique minimal geodesic $\ell$ from $x_t$ to $y_t$, $U(x_t, y_t)$ is the unit tangent vector of $\ell$ at $y_t$, $\xi_t \geq 0$ is a smooth function of $x_t$ to be determined, and

$$\tau := \inf\{t \geq 0 : x_t = y_t\}$$

is the coupling time. Since all terms involved in the equation are regular enough, there exists a unique solution $y_t$. Furthermore, since the additional term containing $1_{\{t<\tau\}}$ vanishes from the coupling time on, one has $x_t = y_t$ for $t \geq \tau$ due to the uniqueness of solutions.

LEMMA 2.3. *Assume that (1.4) and (1.5) hold with $\delta \geq 2\sigma\sqrt{d-1}$. Then there exists a constant $C_3 > 0$ independent of $x, y$ and $T$ such that $x_T = y_T$ holds for $\xi_t := C_3 + 2\sigma\sqrt{d-1}\rho_o(x_t) + \frac{\rho(x,y)}{T}$.*

PROOF. According to Section 2 in [3], we have

$$\begin{aligned}(2.4)\quad d\rho(x_t, y_t) = \ &\{I(x_t, y_t) + \langle \nabla V, \nabla \rho(\cdot, y_t)\rangle(x_t) \\ &+ \langle \nabla V, \nabla \rho(x_t, \cdot)\rangle(y_t) - \xi_t\} \, dt, \qquad t < \tau,\end{aligned}$$

where

$$I_Z(x_t, y_t) = \sum_{i=1}^{d-1} \int_0^{\rho(x_t,y_t)} (|\nabla_U J_i|^2 - \langle R(U, J_i)U, J_i\rangle)(\ell_s) \, ds$$

for $R$ the Riemann curvature tensor, $U$ the unit tangent vector of the minimal geodesic $\ell : [0, \rho(x_t, y_t)] \to M$ from $x_t$ to $y_t$, and $\{J_i\}_{i=1}^{d-1}$ the Jacobi fields along $\ell$, which, together with $U$, consist of an orthonormal basis of the tangent space at $x_t$ and $y_t$ and satisfy

$$J_i(y_t) = P_{x_t,y_t} J_i(x_t), \qquad i = 1, \ldots, d-1.$$

By (1.5) we take a constant $c \geq 0$ such that $\mathrm{Ric} \geq -(c + \sigma^2 \rho_o^2)$. Letting

$$K(x_t, y_t) = \sup_{\ell([0,\rho(x_t,y_t)])} \{c + \sigma^2 \rho_o^2\},$$

we obtain from Wang [21], Theorem 2.14 (see also [7, 8]), that

$$(2.5)\quad I(x_t, y_t) \leq 2\sqrt{K(x_t, y_t)(d-1)} \tanh\left[\frac{\rho(x_t, y_t)}{2}\sqrt{K(x_t, y_t)/(d-1)}\right].$$

Moreover, by (1.4) there exist two constants $r_0, r_1 > 0$ such that $-\text{Hess}_V \geq \delta$ outside $B(o, r_0)$ but $\leq r_1$ on $B(o, r_0)$, where $B(o, r_0)$ is the closed geodesic ball at $o$ with radius $r_0$. Since the length of $\ell$ contained in $B(o, r_0)$ is less than $2r_0$, we conclude that

$$\langle \nabla V, \nabla \rho(\cdot, y_t) \rangle (x_t) + \langle \nabla V, \nabla \rho(x_t, \cdot) \rangle (y_t)$$
$$= \int_0^{\rho(x_t, y_t)} \text{Hess}_V(U, U)(\ell_s)\, ds \leq 2r_0 r_1 - (\rho(x_t, y_t) - 2r_0)^+ \delta$$
$$\leq c_1 - \delta \rho(x_t, y_t)$$

for some constant $c_1 > 0$. Combining this with (2.4), (2.5) and

$$\xi_t = C_3 + 2\sigma\sqrt{d-1}\, \rho_o(x_t) + \frac{\rho(x, y)}{T},$$

we arrive at

$$d\rho(x_t, y_t) \leq \Big\{ 2\sqrt{K(x_t, y_t)(d-1)} + c_1 - \delta\rho(x_t, y_t)$$
$$- C_3 - 2\sigma\sqrt{d-1}\rho_o(x_t) - \frac{\rho(x, y)}{T} \Big\} dt$$

for $t < \tau$. Noting that

$$\sqrt{K(x_t, y_t)} \leq (c + \sigma^2[\rho_o(x_t) + \rho(x_t, y_t)]^2)^{1/2}$$
$$\leq \sqrt{c} + \sigma[\rho_o(x_t) + \rho(x_t, y_t)],$$

and $\delta \geq 2\sigma\sqrt{d-1}$, one has

$$2\sqrt{K(x_t, y_t)(d-1)} - \delta\rho(x_t, y_t) - 2\sigma\sqrt{d-1}\rho_o(x_t) \leq 2\sqrt{c(d-1)}.$$

Thus, when $C_3 \geq c_1 + 2\sqrt{c(d-1)}$ we have

$$d\rho(x_t, y_t) \leq -\frac{\rho(x, y)}{T}\, dt, \qquad t < \tau,$$

so that

$$0 = \rho(x_\tau, y_\tau) \leq \rho(x, y) - \int_0^\tau \frac{\rho(x, y)}{T}\, dt = \frac{T - \tau}{T}\rho(x, y),$$

which implies that $\tau \leq T$ and hence, $x_T = y_T$. □

**3. Harnack inequality and proof of Theorem 1.1.** We first prove the following Harnack inequality using results in Section 2.

PROPOSITION 3.1. *Assume that* (1.4) *and* (1.5) *hold with* $\delta > (1+\sqrt{2})\sigma \times \sqrt{d-1}$. *Then there exist* $C > 0$ *and* $\alpha > 1$ *such that*

$$(3.1) \quad (P_T f(y))^\alpha \leq (P_T f^\alpha(x)) \exp\left[\frac{C}{T}\rho(x,y)^2 + C(T + \rho_o(x)^2)\right]$$

*holds for all* $x, y \in M, T > 0$ *and nonnegative* $f \in C_b(M)$.

PROOF. According to Lemma 2.3, we take

$$\xi_t = C_3 + 2\sigma\sqrt{d-1}\rho_o(x_t) + \frac{\rho(x,y)}{T},$$

such that $\tau \leq T$ and $x_T = y_T$. Obviously, $y_t$ solves the equation

$$d_I y_t = \sqrt{2}\tilde{\Phi}_t \, d\tilde{B}_t + \nabla V(y_t) \, dt$$

for $\tilde{\Phi}_t := P_{x_t,y_t}\Phi_t$ being the horizontal lift of $y_t$, and $\tilde{B}_t$ solving the equation

$$d\tilde{B}_t = dB_t + \frac{1}{\sqrt{2}}\tilde{\Phi}_t^{-1}\xi_t U(x_t, y_t) 1_{\{t < \tau\}} \, dt.$$

By the Girsanov theorem and the fact that $\tau \leq T$, the process $\{\tilde{B}_t : t \in [0,T]\}$ is a $d$-dimensional Brownian motion under the probability measure $R\mathbb{P}$ for

$$R := \exp\left[-\frac{1}{\sqrt{2}}\int_0^\tau \langle P_{x_t,y_t}\Phi_t \, dB_t, \xi_t U(x_t, y_t)\rangle - \frac{1}{4}\int_0^\tau \xi_t^2 \, dt\right].$$

Thus, under this probability measure $\{y_t : t \in [0,T]\}$ is generated by $L$. In particular, $P_T f(y) = \mathbb{E}[f(y_T)R]$. Combining this with the Hölder inequality and noting that $x_T = y_T$, we obtain

$$P_T f(y) = \mathbb{E}[f(y_T)R] = \mathbb{E}[f(x_T)R]$$
$$\leq (P_T f^\alpha(x))^{1/\alpha}(\mathbb{E}R^{\alpha/(\alpha-1)})^{(\alpha-1)/\alpha}.$$

That is,

$$(3.2) \quad (P_T f(y))^\alpha \leq (P_T f^\alpha(x))(\mathbb{E}R^{\alpha/(\alpha-1)})^{\alpha-1}.$$

Since for any continuous exponential integrable martingale $M_t$ and any $\beta, p > 1$, the process $\exp[\beta p M_t - \frac{p^2\beta^2}{2}\langle M\rangle_t]$ is a martingale, by the Hölder inequality one has

$$(3.3) \quad \mathbb{E}e^{\beta M_t - (\beta/2)\langle M\rangle_t} = \mathbb{E}[e^{\beta M_t - (\beta^2 p/2)\langle M\rangle_t} \cdot e^{(\beta(\beta p - 1)/2)\langle M\rangle_t}]$$
$$\leq \mathbb{E}(e^{(\beta p(\beta p - 1)/(2(p-1)))\langle M\rangle_t})^{(p-1)/p}.$$

By taking $\beta = \alpha/(\alpha - 1)$ we obtain

$$(3.4) \quad (\mathbb{E}R^{\alpha/(\alpha-1)})^{\alpha-1}$$
$$\leq \left\{\mathbb{E}\exp\left[\frac{p\alpha(p\alpha - \alpha + 1)}{8(p-1)(\alpha-1)^2}\int_0^T \xi_t^2 \, dt\right]\right\}^{(\alpha-1)(p-1)/p}, \quad p > 1.$$

Since $\delta > (1+\sqrt{2})\sigma\sqrt{d-1}$, we may take $\delta_0 \in ((1+\sqrt{2})\sigma\sqrt{d-1}, \delta)$, small $\varepsilon' > 0$ and large $C_4 > 0$, independent of $T, x$ and $y$, such that

$$\xi_t^2 = \left(C_3 + 2\sigma\sqrt{d-1}\rho_o(x_t) + \frac{\rho(x,y)}{T}\right)^2$$
$$\leq (1-\varepsilon')\left[C_4 + \frac{C_4\rho(x,y)^2}{T^2} + 2(\delta_0 - \sigma\sqrt{d-1})^2\rho_o(x_t)^2\right]$$

holds. Moreover, since

(3.5) $$\lim_{p \downarrow 1} \lim_{\alpha \uparrow \infty} \frac{p\alpha(p\alpha - \alpha + 1)}{8(p-1)(\alpha-1)^2} = \frac{1}{8},$$

there exist $p, \alpha > 1$ such that

$$\frac{p\alpha(p\alpha - \alpha + 1)}{8(p-1)(\alpha-1)^2} \int_0^T \xi_t^2 \, dt$$
$$\leq C_4 T + \frac{C_4\rho(x,y)^2}{T} + \frac{(\delta_0 - \sigma\sqrt{d-1})^2}{4} \int_0^T \rho_o(x_t)^2 \, dt.$$

Combining this with (3.4) and Lemma 2.2, we obtain

$$(\mathbb{E}R^{\alpha/(\alpha-1)})^{\alpha-1} \leq \exp\left[C_5 T + \frac{C_5\rho(x,y)}{T} + C_5\rho_o(x)^2\right], \qquad T > 0, x \in M,$$

for some constant $C_5 > 0$. This completes the proof by (3.2). $\square$

PROOF OF THEOREM 1.1. By Proposition 3.1, let $\alpha > 1$ and $C > 0$ such that (3.1) holds. Since $\delta > \sigma\sqrt{d-1}$, we may take $T > 0$ such that

$$\frac{C}{T} \leq \varepsilon := \frac{1}{8}(\delta - \sigma\sqrt{d-1}).$$

Then for any nonnegative $f \in C_b(M)$ with $\mu(f^\alpha) = 1$, since $\mu$ is $P_T$-invariant, it follows from (3.1) that

$$1 = \int_M P_T f^\alpha(x)\mu(dx) \geq (P_T f(y))^\alpha \int_M e^{-\varepsilon\rho(x,y)^2 - C(1+\rho_o(x)^2)}\mu(dx)$$
$$\geq (P_T f(y))^\alpha \int_{\{\rho_o \leq 1\}} e^{-\varepsilon(1+\rho_o(y))^2 - 2C}\mu(dx)$$
$$\geq \varepsilon'(P_T f(y))^\alpha \exp[-2\varepsilon\rho_o(y)^2], \qquad y \in M,$$

for some constant $\varepsilon' > 0$. Thus,

$$\int_M (P_T f(y))^{2\alpha}\mu(dy) \leq \frac{1}{\varepsilon'}\int_M e^{4\varepsilon\rho_o(y)^2}\mu(dy) < \infty,$$

according to Lemma 2.1. This implies that

$$\|P_T\|_{L^\alpha(\mu) \to L^{2\alpha}(\mu)} < \infty.$$

Therefore, the log-Sobolev inequality (1.3) holds for some constant $C > 0$, due to the uniformly positively improving property of $P_t$ (see [20], proof of Theorem 1.1, and [1]). □

**4. Supercontractivity and ultracontractivity.** Recall that $P_t$ is called supercontractive if $\|P_t\|_{2\to 4} < \infty$ for all $t > 0$ while ultracontractive if $\|P_t\|_{2\to\infty} < \infty$ for all $t > 0$ (see [10]). In the present framework these two properties are stronger than the hypercontractivity: $\|P_t\|_{2\to 4} \leq 1$ for some $t > 0$, which is equivalent to (1.3) due to Gross [14].

PROPOSITION 4.1. *Under* (1.4) *and* (1.5), $P_t$ *is supercontractive if and only if* $\mu(\exp[\lambda\rho_o^2]) < \infty$ *for all* $\lambda > 0$, *while it is ultracontractive if and only if* $\|P_t \exp[\lambda\rho_o^2]\|_\infty < \infty$ *for all* $t, \lambda > 0$.

PROOF. The proof is similar to that of [18], Theorem 2.3. Let $f \in L^2(\mu)$ with $\mu(f^2) = 1$. By (3.1) for $\alpha = 2$ and noting that $\mu$ is $P_t$-invariant, we obtain

$$1 \geq (P_T f(y))^2 \int_M \exp\left[-\frac{C}{T}\rho(x,y)^2 - C(T + \rho_o(x)^2)\right]\mu(dx)$$

$$\geq (P_T f(y))^2 \exp\left[-\frac{2C}{T}(\rho_o(y)^2 + 1) - C(T+1)\right]\mu(B(o,1)).$$

Hence, for any $T > 0$ there exists a constant $\lambda_T > 0$ such that

(4.1) $\qquad |P_T f| \leq \exp[\lambda_T(1 + \rho_o^2)], \qquad T > 0, \mu(f^2) = 1.$

(1) If $\mu(e^{\lambda\rho_o^2}) < \infty$ for any $\lambda > 0$, (4.1) yields that

$$\|P_T\|_{2\to 4}^4 \leq \mu(e^{4\lambda_T(1+\rho_o^2)}) < \infty, \qquad T > 0.$$

Conversely, if $P_t$ is supercontractive then the super log-Sobolev inequality (cf. [10])

$$\mu(f^2 \log f^2) \leq r\mu(|\nabla f|^2) + \beta(r), \qquad r > 0, \mu(f^2) = 1,$$

holds for some $\beta:(0,\infty) \to (0,\infty)$. By [2] (see also [17, 18]), this inequality implies $\mu(e^{\lambda\rho_o^2}) < \infty$ for all $\lambda > 0$.

(2) By (4.1) and the semigroup property,

$$\|P_T\|_{2\to\infty} \leq \|P_{T/2} e^{\lambda_{T/2}(1+\rho_o^2)}\|_\infty < \infty, \qquad T > 0,$$

provided $\|P_t e^{\lambda\rho_o^2}\|_\infty < \infty$ for any $t, \lambda > 0$. Conversely, since the ultracontractivity is stronger than the supercontractivity, it implies that $e^{\lambda\rho_o^2} \in L^2(\mu)$ for any $\lambda > 0$ as explained above. Therefore,

$$\|P_t e^{\lambda\rho_o^2}\|_\infty \leq \|P_t\|_{2\to\infty} \|e^{\lambda\rho_o^2}\|_2 < \infty, \qquad \lambda > 0.$$

Then the proof is completed. □

To derive explicit conditions for the supercontractivity and ultracontractivity, we consider the following stronger version of (1.4):

(4.2) $\quad -\mathrm{Hess}_V \geq \Phi \circ \rho_o \quad$ holds outside a compact subset of $M$

for a positive increasing function $\Phi$ with $\Phi(r) \uparrow \infty$ as $r \uparrow \infty$. We then aim to search for reasonable conditions on positive increasing function $\Psi$ such that

(4.3) $$\mathrm{Ric} \geq -\Psi \circ \rho_o$$

implies the supercontractivity and/or ultracontractivity.

THEOREM 4.2. *If (4.3) and (4.2) hold for some increasing positive functions $\Phi$ and $\Psi$ such that*

(4.4) $$\lim_{r\to\infty} \Phi(r) = \lim_{r\to\infty} \frac{(\int_0^r \Phi(s)\,ds)^2}{\Phi(r)} = \infty,$$

(4.5) $$\sqrt{\Psi(r+t)(d-1)} \leq \theta \int_0^r \Phi(s)\,ds + \frac{1}{2}\int_0^{t/2} \Phi(s)\,ds + C, \qquad r,t \geq 0,$$

*for some constants $\theta \in (0, 1/(1+\sqrt{2}))$ and $C > 0$. Then $P_t$ is supercontractive. Furthermore, if*

(4.6) $$\int_1^\infty \frac{ds}{\sqrt{s}\int_0^{\sqrt{r}} \Phi(u)\,du} < \infty,$$

*then $P_t$ is ultracontractive. More precisely, for*

$$\Gamma_1(r) := \frac{1}{\sqrt{r}}\int_0^{\sqrt{r}} \Phi(s)\,ds, \qquad \Gamma_2(r) := \int_r^\infty \frac{ds}{\sqrt{s}\int_0^{\sqrt{s}} \Phi(u)\,du}, \qquad r > 0,$$

*(4.6) implies*

(4.7) $\|P_t\|_{2\to\infty} \leq \exp\left[c + \frac{c}{t}(1 + \Gamma_1^{-1}(c/t) + \Gamma_2^{-1}(t/c))\right] < \infty, \qquad t > 0,$

*for some constant $c > 0$ and*

$$\Gamma_1^{-1}(s) := \inf\{t \geq 0 : \Gamma_1(t) \geq s\}, \qquad s \geq 0.$$

PROOF. (a) Replacing $c+\rho_o^2$ by $\Psi\circ\rho_o$ and noting that $\text{Hess}_V \leq -\Phi\circ\rho_o$ for large $\rho_o$, the proof of Lemma 2.1 implies

$$(4.8) \qquad L\rho_o^2 \leq c_1(1+\rho_o) - 2\rho_o\left(\int_0^{\rho_o}\Phi(s)\,ds - \sqrt{\Psi\circ\rho_o(d-1)}\right)$$

for some constant $c_1 > 0$. Combining this with (4.5) and noting that $\frac{1}{\rho_o} \times \int_0^{\rho_o}\Phi(s)\,ds \to \infty$ as $\rho_o \to \infty$, we conclude that for any $\lambda > 0$,

$$\begin{aligned}(4.9) \qquad Le^{\lambda\rho_o^2} &\leq C - \frac{2\lambda\rho_o\sqrt{2}}{1+\sqrt{2}}e^{\lambda\rho_o^2}\int_0^{\rho_o}\Phi(s)\,ds + 4\lambda^2\rho_o^2 e^{\lambda\rho_o^2} \\ &\leq C + C(\lambda) - \lambda\rho_o e^{\lambda\rho_o^2}\int_0^{\rho_o}\Phi(s)\,ds,\end{aligned}$$

where $C > 0$ is a universal constant and

$$\begin{aligned}(4.10) \qquad C(\lambda) &:= \sup_{r>0} re^{\lambda r^2}\left\{4\lambda^2 r - \frac{\lambda}{(1+\sqrt{2})^2}\int_0^r\Phi(s)\,ds\right\} \\ &= \sup_{r^2\leq\Gamma_1^{-1}(4(1+\sqrt{2})^2\lambda)} re^{\lambda r^2}\left\{4\lambda^2 r - \frac{\lambda}{(1+\sqrt{2})^2}\int_0^r\Phi(s)\,ds\right\} \\ &\leq 4\lambda^2\Gamma_1^{-1}(4(1+\sqrt{2})^2\lambda)\exp[\lambda\Gamma_1^{-1}(4(1+\sqrt{2})^2\lambda)] \\ &\leq \exp[4\lambda + 2\lambda\Gamma_1^{-1}(4(1+\sqrt{2})^2\lambda)] < \infty.\end{aligned}$$

Therefore, (1.1) holds and

$$(4.11) \qquad \mu(e^{\lambda\rho_o^2}) < \infty, \qquad \lambda > 0.$$

(b) By (4.5), (4.8) and Kendall's Itô formula [16] as in the proof of Lemma 2.2, we have

$$d\rho_o^2(x_t) \leq 2\sqrt{2}\rho_o(x_t)\,db_t + \left(C_1 - \frac{2\sqrt{2}\rho_o(x_t)(1+\varepsilon)}{1+\sqrt{2}}\int_0^{\rho_o(x_t)}\Phi(s)\,ds\right)dt$$

for some constants $\varepsilon, C_1 > 0$, where $x_t$ and $b_t$ are in the proof of Lemma 2.2. Let

$$(4.12) \qquad \varphi(r) = \int_0^r \frac{ds}{\sqrt{s}}\int_0^{\sqrt{s}}\Phi(u)\,du, \qquad r \geq 0.$$

We arrive at

$$\begin{aligned}d\varphi\circ\rho_o^2(x_t) &\leq 2\sqrt{2}\rho_o(x_t)\varphi'\circ\rho_o^2(x_t)\,db_t + 4\rho_o^2(x_t)\varphi''\circ\rho_o^2(x_t)\,dt \\ &\quad + \varphi'\circ\rho_o^2(x_t)\left(C_1 - \frac{2\sqrt{2}\rho_o(x_t)(1+\varepsilon)}{1+\sqrt{2}}\int_0^{\rho_o(x_t)}\Phi(s)\,ds\right)dt.\end{aligned}$$

From (4.4) we see that

$$\frac{\rho_o \varphi'' \circ \rho_o^2}{\varphi' \circ \rho_o^2 \int_0^{\rho_o} \Phi(s)\, ds} \leq \frac{\Phi \circ \rho_o}{2(\int_0^{\rho_o} \Phi(s)\, ds)^2},$$

which goes to zero as $\rho_o \to \infty$. Then there exists a constant $C_2 > C_1$ such that

$$d\varphi \circ \rho_o^2(x_t) \leq 2\sqrt{2}\left(\int_0^{\rho_o(x_t)} \Phi(s)\, ds\right) db_t$$
$$+ C_2\, dt - \frac{2\sqrt{2}}{1+\sqrt{2}}\left(\int_0^{\rho_o(x_t)} \Phi(s)\, ds\right)^2 dt.$$

This implies that for any $\lambda > 0$,

$$\mathbb{E}\exp\left[\frac{2\sqrt{2}\,\lambda}{1+\sqrt{2}}\int_0^T \left(\int_0^{\rho_o(x_t)} \Phi(s)\, ds\right)^2 dt\right]$$
$$\leq e^{C_2 \lambda T + \lambda \varphi \circ \rho_o^2(x)} \mathbb{E}\exp\left[2\sqrt{2}\lambda \int_0^T \left(\int_0^{\rho_o(x_t)} \Phi(s)\, ds\right) db_t\right]$$
$$\leq e^{C_2 \lambda T + \lambda \varphi \circ \rho_o^2(x)} \left(\mathbb{E}\exp\left[16\lambda^2 \int_0^T \left(\int_0^{\rho_o(x_t)} \Phi(s)\, ds\right)^2 dt\right]\right)^{1/2}.$$

Taking

$$\lambda = \frac{\sqrt{2}}{8(1+\sqrt{2})},$$

we arrive at

$$\mathbb{E}\exp\left[\frac{1}{2(1+\sqrt{2})^2}\int_0^T \left(\int_0^{\rho_o(x_t)} \Phi(s)\, ds\right)^2 dt\right]$$

(4.13)
$$\leq e^{2C_2 T + \varphi \circ \rho_o^2(x)\sqrt{2}/8(1+\sqrt{2})}.$$

(c) Let $\gamma\colon [0, \rho(x_t, y_t)] \to M$ be the minimal geodesic from $x_t$ to $y_t$, and $U$ its tangent unit vector. By (4.2), there exists a constant $C_3 > 0$ such that

(4.14)
$$\langle \nabla V, \nabla \rho(\cdot, y_t)\rangle(x_t) + \langle \nabla V, \nabla \rho(x_t, \cdot)\rangle(y_t)$$
$$= \int_0^{\rho(x_t, y_t)} \mathrm{Hess}_V(U_s, U_s)\, ds \leq C_3 - \int_0^{\rho(x_t, y_t)/2} \Phi(s)\, ds.$$

To understand the last inequality, we assume, for instance, that $\rho_o(x_t) \geq \rho_o(y_t)$ so that by the triangle inequality,

$$\rho_o(\gamma_s) \geq \rho_o(x_t) - s \geq \rho(x_t, y_t)/2 - s, \qquad s \in [0, \rho(x_t, y_t)/2].$$

For the coupling constructed in Section 3, one concludes from (4.14) and the proof of Lemma 2.3 that

$$d\rho(x_t, y_t) \leq \left\{ 2\sqrt{K(x_t, y_t)(d-1)} + C_4 \right.$$
(4.15)
$$\left. - \int_0^{\rho(x_t,y_t)/2} \Phi(s)\, ds - \xi_t \right\} dt, \qquad t < \tau,$$

holds for some constant $C_4 > 0$, where

$$K(x_t, y_t) := \sup_{\ell([0, \rho(x_t, y_t)])} \Psi \circ \rho_o \leq \Psi(\rho_o(x_t) + \rho(x_t, y_t)),$$

and $\ell$ is the minimal geodesic from $x_t$ to $y_t$. Combining (4.5) and (4.15), we obtain

$$d\rho(x_t, y_t) \leq \left\{ C_4 + 2\theta \int_0^{\rho_o(x_t)} \Phi(s)\, ds - \xi_t \right\} dt, \qquad t < \tau.$$

So, taking

$$\xi_t = C_4 + 2\theta \int_0^{\rho_o(x_t)} \Phi(s)\, ds + \frac{\rho(x, y)}{T},$$

we arrive at

$$d\rho(x_t, y_t) \leq -\frac{\rho(x, y)}{T}\, dt, \qquad t < \tau.$$

This implies $\tau \leq T$, and hence $x_T = y_T$ a.s.

Combining (4.5) with (3.4) and (3.5) we conclude that for the present choice of $\xi_t$ there exist $\alpha, p, C_5 > 1$ such that

$$(\mathbb{E}R^{\alpha/(\alpha-1)})^{p/(p-1)} \leq \mathbb{E}\exp\left[\frac{1}{2(1+\sqrt{2})^2} \int_0^T \left(\int_0^{\rho(x_t)} \Phi(s)\, ds\right)^2 dt \right.$$
$$\left. + C_5 T + \frac{C_5}{T}\rho(x,y)^2 \right].$$

Combining this with (4.13) and (3.2) we obtain

(4.16) $\quad (P_T f(y))^\alpha \leq (P_T f^\alpha(x)) \exp\left[ CT + \frac{C}{T}\rho(x,y)^2 + C\varphi \circ \rho^2(x) \right]$

holds for some $\alpha, C > 1$, any positive $f \in C_b(M)$ and all $x, y \in M, T > 0$.

(d) For any positive $f \in C_b(M)$ with $\mu(f^\alpha) = 1$, (4.16) implies that

$$(P_T f(y))^\alpha \int_{B(o,1)} \exp\left[-CT - \frac{C}{T}\rho(x,y)^2 - C\varphi^2(x)\right] \mu(dx) \leq 1.$$

Therefore, there exists a constant $C' > 0$ such that

$$(4.17) \quad (P_T f(y))^\alpha \leq \exp\left[C'(1+T) + \frac{C'}{T}\rho(y)^2\right], \qquad y \in M, T > 0.$$

Combining this with (4.11) we obtain

$$\|P_T\|_{\alpha \to p\alpha} < \infty, \qquad T > 0, p > 1.$$

This is equivalent to the supercontactivity by the Riesz–Thorin interpolation theorem and $\|P_t\|_{1\to 1} = 1$. Thus, the first assertion holds.

(e) To prove (4.7), it suffices to consider $t \in (0, 1]$ since $\|P_t\|_{2\to\infty}$ is decreasing in $t > 0$. So, below we assume that $T \leq 1$. By (4.17) and the fact that $(P_{2T} f)^\alpha \leq P_T (P_T f)^\alpha$, we have

$$(4.18) \qquad \|P_{2T}\|_{\alpha \to \infty} \leq \|P_T e^{2C'\rho_o^2/T}\|_\infty e^{C'(1+T)}, \qquad T > 0.$$

Therefore, by the Riesz–Thorin interpolation theorem and $\|P_t\|_{1\to 1} = 1$, for the ultracontractivity it suffices to show that

$$(4.19) \qquad\qquad \|P_T e^{\lambda \rho_o^2}\|_\infty < \infty, \qquad \lambda, T > 0.$$

Since $\Phi$ is increasing, it is easy to check that

$$\eta(r) := \sqrt{r} \int_0^{\sqrt{r}} \Phi(s)\, ds, \qquad r \geq 0,$$

is convex, and so is $s \mapsto s\eta(\frac{\log s}{\lambda})$ for $\lambda > 0$. Thus, it follows from (4.9) and the Jensen inequality that

$$h_{\lambda,x}(t) := \mathbb{E} e^{\lambda \rho_o^2(x_t)} < \infty, \qquad x_0 = x \in M, \lambda, t > 0,$$

and

$$\frac{d^+}{dt} h_{\lambda,x}(t) \leq C + C(\lambda) - \lambda h_{\lambda,x}(t) \eta(\lambda^{-1} \log h_{\lambda,x}(t)), \qquad t > 0.$$

This implies (4.19), provided (4.6) holds. This can be done by considering the following two situations:

(1) Since $h_{\lambda,x}(t)$ is decreasing provided $\lambda h_{\lambda,x}(t)\eta(\lambda^{-1} \log h_{\lambda,x}(t)) > C + C(\lambda)$, if

$$\lambda h_{\lambda,x}(0)\eta(\lambda^{-1}\log h_{\lambda,x}(0)) \leq 2C + 2C(\lambda),$$

then

$$h_{\lambda,x}(t) \leq \sup\{r \geq 1 : \lambda r \eta(\lambda^{-1}\log r) \leq 2C + 2C(\lambda)\} \leq \frac{1}{\lambda}(2C + 2C(\lambda)) + C''$$

for some constant $C'' > 0$.

(2) If $\lambda h_{\lambda,x}(0)\eta(\lambda^{-1}\log h_{\lambda,x}(0)) > 2C + 2C(\lambda)$, then $h_{\lambda,x}(t)$ is decreasing in $t$ up to

$$t_\lambda := \inf\{t \geq 0 : \lambda h_{\lambda,x}(t)\eta(\lambda^{-1}\log h_{\lambda,x}(t)) \leq 2C + 2C(\lambda)\}.$$

Indeed,

$$\frac{d^+}{dt}h_{\lambda,x}(t) \leq -\frac{\lambda}{2}h_{\lambda,x}(t)\eta(\lambda^{-1}\log h_{\lambda,x}(t)), \qquad t \leq t_\lambda.$$

Thus,

$$\int_{h_{\lambda,x}(T\wedge t_\lambda)}^{\infty} \frac{dr}{r\eta(\lambda^{-1}\log r)} \geq \frac{\lambda}{2}(T \wedge t_\lambda).$$

This is equivalent to

$$\Gamma_2(\lambda^{-1}\log h_{\lambda,x}(T \wedge t_\lambda)) \geq \tfrac{1}{2}(T \wedge t_\lambda).$$

Hence,

$$h_{\lambda,x}(T \wedge t_\lambda) \leq \exp[\lambda \Gamma_2^{-1}(\tfrac{1}{2}(T \wedge t_\lambda))].$$

Since it is reduced to case (1) if $T > t_\lambda$ by regarding $t_\lambda$ as the initial time, in conclusion we have

$$\sup_{x \in M} h_{\lambda,x}(T) \leq \max\Big\{\exp[\lambda\Gamma_2^{-1}(T/2)],\ C'' + \frac{1}{\lambda}(2C + 2C(\lambda))\Big\}.$$

Therefore, (4.7) follows from (4.18), (4.10) with $\lambda = 2C'/T$, and the Riesz interpolation theorem. $\square$

Finally, we note that a simple example for conditions in Theorem 4.2 to hold is

$$\Phi(s) = s^{\alpha-1}, \qquad \Psi(s) = \varepsilon s^{2\alpha}$$

for $\alpha > 1$ and small enough $\varepsilon > 0$. In this case $P_t$ is ultracontractive with

$$\|P_t\|_{2\to\infty} \leq \exp[c(1 + t^{-(\alpha+1)/(\alpha-1)})], \qquad t > 0,$$

for some $c > 0$.

**Acknowledgments.** The author would like to thank the referees for their careful reading and valuable comments on an earlier version of the paper.

School of Mathematics  
Beijing Normal University  
Beijing 100875  
China  
and  
Swansea University  
Singleton Park  
Swansea, SA2 8PP  
United Kingdom  
E-mail: wangfy@bnu.edu.cn